\documentclass[oneside]{amsart}

\usepackage{control/packages}
\usepackage{control/docaesthetic}
\usepackage{control/environments}
\usepackage{control/macros}
\usepackage{control/paperdata}

\usepackage{subfiles}

\usepackage[style=numeric]{biblatex}
\begin{document}

\maketitle

%%%%
%%%
%%
%

\begin{abstract}
We show that pointlike sets are decidable for the pseudovariety of finite semigroups whose idempotent-generated subsemigroup is ${\mathcal R}$-trivial. 
Notably, our proof is constructive: we provide an explicit relational morphism which computes the ${\mathbf{ER}}$-pointlike subsets of a given finite semigroup.
\end{abstract}

%
%%
%%%
%%%%

\tableofcontents

%%%%
%%%
%%
%

\section{Introduction}

Let $\pvV$ be a pseudovariety.
A non-empty subset $X$ of a finite semigroup $S$ is said to be \textit{$\pvV$-pointlike} if for any relational morphism $\rho : S \rlm V$ with $V \in \pvV$ there exists $v \in V$ for which $X \subseteq (v)\rho^\inv$.
If there is an algorithm which produces the $\pvV$-pointlike subsets of any finite semigroup given as input, then $\pvV$ is said to have \textit{decidable pointlikes}.
For background on pointlike sets, see \cite{STEINBERG-PL-SURVEY, GENTHEORYPL}.

The main result of this paper (Theorem \ref{thm:maintheorem}) is that the pseudovariety
    \[
        \pvER = \setst{ S \in \FSGP }{\langle E(S) \rangle \in \pvR}
    \]
has decidable pointlikes, where $\pvR$ is the pseudovariety of $\GR$-trivial semigroups.
Moreover, this result is proven \textit{constructively}; that is, for any finite semigroup $S$, we can explicitly construct a relational morphism from $S$ to a member of $\pvER$ which computes the $\pvER$-pointlike subsets of $S$.  

\subsection{Organization of paper}
Section \ref{section:prelim} covers various preliminary notions and establishes basic notation.
Basic facts about $\pvER$---particularly relating to group kernels and the type-II partition---are covered in Section \ref{section:type2}.

Section \ref{section:framework} provides an overview of various key notions from the authors' framework for pointlike sets---which was developed in \cite{GENTHEORYPL}---and establishes the language used in the rest of the paper.
This framework is particularized in Section \ref{section:lowerbound} to define a candidate definition for $\pvER$-pointlikes and to prove that said candidate provides a lower bound.

Section \ref{section:stable} establishes a number of preliminary results which are necessary for the construction defined in Section \ref{section:automaton}.
The primary data of said construction---and the central ``conceptual contribution'' of the paper---is an automaton whose transition semigroup is shown in Section \ref{section:upperbound} to belong to $\pvER$.\footnotemark
\footnotetext{This automaton---which the first author has nicknamed \textit{PFFL} (\underline{P}ermute \underline{F}irst, \underline{F}all \underline{L}ater)---could be seen as an evolution of a construction used in \cite{CPL-LOWERBOUNDS}.}
This establishes by way of the previous two sections that our lower bound for $\pvER$-pointlikes is an upper bound as well, which in turn establishes the main result (Theorem \ref{thm:maintheorem}).

%
%%
%%%
%%%%

%%%%
%%%
%%
%

\section{Preliminaries} \label{section:prelim}

\subsection{}
Familiarity with finite semigroup theory and basic category theory is assumed; for reference, the reader is directed to \cite{QTHEORY} (for finite semigroup theory) and \cite{CTCONTEXT} (for category theory).

\begin{notation}
Let $S$ be a finite semigroup.
\begin{itemize}
\item Write $S^I$ for the semigroup obtained by adjoining a new element $I$ to $S$ and defining $xI = Ix = x$ for all $x \in S$.
\item Let $E(S)$ denote the set of idempotents of $S$.
\item Given $x \in S$, let $x^\omega$ denote the unique idempotent generated by $x$.
\item Green's equivalence relations are denoted by $\GR$, $\GL$, $\GH$, and $\GJ$; 
moreover, the various Green's equivalence classes of $x \in S$ are denoted by $\RCL{x}$, $\LCL{x}$, $\HCL{x}$, and $\JCL{x}$, respectively.
\end{itemize}
\end{notation}

\subsection{Partial transformation semigroups}
A (finite) \textbf{partial transformation semigroup}, which we will abbreviate as PTS, is a pair $(Q , S)$ consisting of a finite set $Q$ and a finite semigroup $S$ which acts on the right of $Q$ by partial functions.
If $q \ast s$ is undefined for some $q \in Q$ and some $s \in S$, we will write $q \ast s = \varnothing$.

A PTS morphism $(\zeta , \varphi) : (Q , S) \rightarrow (P , T)$ is given by a pair 
    \[
        \zeta : Q \arw P
            \quad \text{and} \quad
        \varphi : S \arw T
    \]
where $\zeta$ is a set function and $\varphi$ is a morphism of semigroups such that
    \[
        (q)\zeta \ast (s)\varphi = (q \ast s)\zeta
    \]
for all $q \in Q$ and all $s \in S$ for which $q \ast s \neq \varnothing$.

A \textbf{congruence} on a PTS $(Q , S)$ is an equivalence relation $\equiv$ on $Q$ such that
    \[
        q_1 \equiv q_2
            \quad \Longrightarrow \quad
        q_1 \ast s \equiv q_2 \ast s
    \]
for all $q_1, q_2 \in Q$ and all $s \in S$ for which both $q_1 \ast s$ and $q_2 \ast s$ are defined.
Such a congruence $\equiv$ on $(Q , S)$ induces a quotient PTS
    \[
        (Q , S) / \equiv 
            \;\;=\;\;
        ( Q / \equiv , S )
    \]
wherein the action is given by
    \[
        [q]_\equiv \ast s
            =
        \begin{cases}
            [ q \ast s ]_\equiv,     &   \text{if $q' \ast s \neq \varnothing$ for some $q' \equiv q$;} \\
            \varnothing,             &   \text{otherwise;}
        \end{cases}
    \]
for each $q \in Q$ and each $s \in S$.

A PTS $(Q , S)$ is said to be \textbf{injective} if $S$ acts on $Q$ by partial injections.
Moreover, a congruence $\equiv$ on a PTS $(Q , S)$ is said to be \textit{injective} if $(Q , S) / \equiv$ is injective.

A (finite) \textbf{transformation semigroup} is a PTS $(Q , S)$ wherein the action of $S$ on $Q$ is by fully defined functions on $Q$.
If a transformation semigroup $(Q , S)$ satisfies the condition that 
    \[
        \text{$q \ast s_1 = q \ast s_2$ for all $q \in Q$} 
            \quad \Longrightarrow \quad
        s_1 = s_2
    \]
for all $s_1 , s_2 \in S$, then $(Q , S)$ is said to be \textbf{faithful}.

\subsection{Activators}
Let $J$ be a $\GJ$-class of a finite semigroup $S$.
The set 
    \[
    \setst{a \in S^I}{Ja \cap J \neq \varnothing}
    \]
is a union of $\GJ$-classes of $S^I$ which contains a unique $\lgJ$-minimal $\GJ$-class.
Said minimal $\GJ$-class---which is always regular---is called the \textbf{right activator} of $J$ and will be denoted by $\RACT(J)$.
Note that $J$ is regular if and only if $\RACT(J) = J$.

\begin{lem} \label{prelim:lem:activator}
Let $J$ be a $\GJ$-class of a finite semigroup $S$.
For each $x \in J$, there exists an element $t \in \RACT(J)$ such that
\begin{enumerate}
    \item $xt = x$;
    \item $xs \slgR x$ if and only if $ts \slgR t$ for any $s \in S$; and
    \item left multiplication by $x$ defines a surjective function from the $\GR$-class of $t$ onto the $\GR$-class of $x$.
\end{enumerate}
In fact, $t$ may be chosen to be an idempotent.
\end{lem}

\begin{proof}
See \cite[Lemma~2.9]{CPL-LOWERBOUNDS}.
\end{proof}

\subsection{}
Given $x \in S$, let $\cF_x$ denote the set of elements in $\RACT( \JCL{x} )$ which satisfy the claims of Lemma \ref{prelim:lem:activator}.
Note that if $x$ is regular then $\cF_x$ is the set of idempotents which are $\GL$-equivalent to $x$.

\subsection{Relational morphisms}
A \textbf{relational morphism} $\rho : S \rlm T$ is an equivalence class of spans in the category of finite semigroups of the form
    \begin{equation*}
        \begin{tikzcd}
            \cdot 
            \arrow[r] 
            \arrow[d, two heads] 
        & 
            T 
        \\
            S 
        &
        \end{tikzcd}
    \end{equation*}
where the map to $S$ is a regular epimorphism,\footnotemark and where two such spans are equivalent if the natural maps from each apex to $S \times T$ have the same image.
\footnotetext{Here and throughout, ``regular epimorphism'' means ``surjective homomorphism''.}

\subsection{Pseudovarieties}
A \textbf{pseudovariety} is a class of finite semigroups which is closed under taking subsemigroups, quotients, and finite products of its members.

\subsection{Power semigroups}   \label{prelim:power}
Given a finite semigroup $S$, let $\Po(S)$ denote the semigroup of non-empty subsets of $S$ under the inherited operation given by
    \[
        X \cdot Y = \setst{xy}{x \in X, y \in Y}
    \]
for all non-empty subsets $X$ and $Y$ of $S$; also, let $\sing(S)$ denote the subsemigroup of $\Po(S)$ consisting of the singletons.

A morphism $\varphi : S \rightarrow T$ extends to a morphism
    \[
        \ext{\varphi} : \Po(S) \arw \Po(T)
            \quad \text{given by} \quad
        (X)\ext{\varphi} = \setst{(x)\varphi}{x \in X}.
    \]
Equipping the object map $\Po$ with this action on morphisms yields a functor
    \[
        \Po : \FSGP \arw \FSGP
    \]
which creates monomorphisms, regular epimorphisms, and isomorphisms.

\subsection{Pointlikes}
Given a finite semigroup $S$ and a pseudovariety $\pvV$, a non-empty subset $X \subseteq S$ is said to be \textbf{$\pvV$-pointlike} if for any relational morphism of the form $\rho : S \rlm V$ with $V \in \pvV$ there exists some element $v \in V$ for which $X \subseteq (v)\rho^\inv$.

The set of $\pvV$-pointlike subsets of $S$ is denoted by $\PV(S)$, and is in fact a subsemigroup of $\Po(S)$ which contains $\sing(S)$ and which is closed under taking non-empty subsets of its members.
Equipping this object map with the action on morphisms sending $\varphi : S \rightarrow T$ to the evident restriction of the extension described in \ref{prelim:power} yields a subfunctor
    \[
        \PV : \FSGP \arw \FSGP
    \]
of $\Po$ with the property that a finite semigroup $S$ belongs to $\pvV$ if and only if $\PV(S) = \sing(S)$.\footnotemark
\footnotetext{The notation for $\Po$ is due to it being the pointlikes functor for the trivial pseudovariety $\pvtriv$.}
Pointlike functors also create monomorphisms, regular epimorphisms, and isomorphisms.

%
%%
%%%
%%%%

%%%%
%%%
%%
%

\section{Group kernels and the type-II partition}   \label{section:type2}

\subsection{Group kernels}
Let $S$ be a finite semigroup.
Recall that the \textbf{group kernel} of $S$ is the subsemigroup $\KG(S)$ consisting of those elements which are always contained in the inverse image of the identity under any relational morphism from $S$ to a finite group.
That is, $x \in \KG(S)$ if and only if $x \in (1_G)\rho^\inv$ for any relational morphism $\rho : S \rlm G$ with $G \in \pvG$.

It is well-known that $\KG(S)$ is the smallest subsemigroup of $S$ for which $E(S) \subseteq \KG(S)$ and such that if $s \in \KG(S)$ and $x, y \in S$ with $xyx=x$, then $xsy, ysx \in \KG(S)$ as well.

\begin{lem}     \label{type2:lem:kgepi}
The object map $\KG$ is a functor when equipped with the evident restriction action on morphisms.
Moreover, $\KG$ preserves regular epimorphisms.
\end{lem}

\begin{proof}
See \cite[Proposition~4.12.6]{QTHEORY}.
\end{proof}

\begin{defn}
The (right-sided) \textbf{type-II partition} on $S$ is defined by
    \[
        x \equiv_\II y
            \quad \Longleftrightarrow \quad
        \text{$xa=y$ and $yb=x$ for some $a,b \in \KG(S)^I$}
    \]
for all $x , y \in S$.
Let $\BII{x}$ denote the type-II equivalence class of $x \in S$.
\end{defn}

\subsection{}
It is easy to see that $\equiv_\II$ is contained in $\GR$.
Given an $\GR$-class $R$ of $S$, consider the PTS $(R, S)$.
The relation $\equiv_\II$ is a congruence on $(R, S)$---that is, given $x, y \in R$ with $x \equiv_\II y$, then $xs \equiv_\II ys$ for any $s \in S$ such that both $xs$ and $ys$ are defined.
Let $(R/\II , S)$ denote the quotient of $(R , S)$ by $\equiv_\II$.
Crucially, $(R/\II , S)$ is an injective PTS---in fact, $\equiv_\II$ is the minimal injective congruence on $(R, S)$.

\begin{lem} \label{type2:lem:mininj}
Let $R$ be an $\GR$-class of a finite semigroup $S$.
Then $\equiv_\II$ is the minimal injective congruence on $(R , S)$.
\end{lem}

\begin{proof}
See \cite{TYPEIIREDUX}.
\end{proof}

\begin{lem}     \label{type2:lem:type2KG}
If $R$ is an $\GR$-class of $S$, then
\begin{enumerate}
    \item $\KG(S) \cap R$ is either empty or a $\II$-class of $R$, and
    \item any $a \in \KG(S)$ acts as a partial identity in $(R/\II , S)$.
\end{enumerate}
\end{lem}

\begin{proof}
Straightforward.
\end{proof}

\begin{lem}     \label{type2:lem:actII}
Let $x \in S$ and let $t \in \cF_x$. Then
\begin{enumerate}
\item $x \cdot \BII{t} = \BII{x}$;
\item if $s \gR t$ then $x \cdot \BII{s} = \BII{xs}$; and
\item there is a PTS morphism
    \[
        ( x \cdot (-) , \id{S} ) : (\RCL{t} / \II , S) \longsurj (\RCL{x}/\II , S),
    \]
sending $\BII{s} \in \RCL{t} / \II$ to $\BII{xs}$ and acting as identity on $S$.
\end{enumerate}
\end{lem}

\begin{proof}
Straightforward.
\end{proof}

\begin{lem}     \label{type2:lem:ERchar}
A finite semigroup $S$ belongs to $\pvER$ if and only if $(R , S)$ is an injective PTS for every $\GR$-class $R$ of $S$.
Moreover, 
    \[
        \pvER = \pvR \ast \pvG = \pvR \malcev \pvG.
    \]
\end{lem}

\begin{proof}
See \cite[Theorem~4.8.3]{QTHEORY}.
\end{proof}

%
%%
%%%
%%%%

%%%%
%%%
%%
%

\section{General theory of pointlike sets}  \label{section:framework}

\subsection{}
In this section we will briefly cover key aspects of the authors' ``general theory of pointlike sets'' which provides the framework for our work here.
The treatment here is incomplete and proofs are omitted; for further details, see \cite{GENTHEORYPL}.

\subsection{Semigroup complexes}
Let $S$ be a finite semigroup.
An \textbf{$S$-complex} is a subsemigroup $\cK \subseteq \Po(S)$ which 
\begin{enumerate}
    \item contains $\sing(S)$, and which
    \item is closed under taking non-empty subsets of its members, i.e., if $X \in \cK$ then any $Y \in \Po(S)$ for which $Y \subseteq X$ also belongs to $\cK$.
\end{enumerate}
The set of $S$-complexes---denoted by $\Com{S}$---is a complete lattice wherein the order is inclusion, the top and bottom are $\Po(S)$ and $\sing(S)$ respectively, the meet is intersection, and the join is given by
    \[
        \cK_1 \vee \cK_2 
            =
        \setst{ X \in \Po(S) }{ \text{$X \subseteq \Tilde{X}$ for some $\Tilde{X} \in \langle \cK_1 \cup \cK_2 \rangle$}}
    \]
for any $\cK_1 , \cK_2 \in \Com{S}$.

\begin{defn}
A \textbf{modulus} $\Lambda$ is a rule which assigns to each finite semigroup $S$ a set $\Lambda_S \subseteq \Po(S)$ in a manner which satisfies the following axioms.
\begin{enumerate}
    \item If $\varphi : S \rightarrow T$ is a morphism, then for any $X \in \Lambda_S$ there exists some $\widetilde{X} \in \Lambda_T$ such that $(X)\ext{\varphi} \subseteq \widetilde{X}$.
    \item If $\varphi : S \surj T$ is a regular epimorphism, then for any $Y \in \Lambda_T$ there exists some $\widetilde{Y} \in \Lambda_S$ such that $(\widetilde{Y})\ext{\varphi} = Y$.
\end{enumerate}
When defining moduli, we will generally write
    \[
        \Lambda = \modrl{S}{\Lambda_S},
    \]
to mean ``$\Lambda$ is the rule which assigns $\Lambda_S$ to a given finite semigroup $S$''.
\end{defn}

\subsection{Constructing lower bounds for pointlikes}   \label{framework:lowerbound}
Given a modulus $\Lambda$, the \textbf{$\Lambda$-construct} of a finite semigroup $S$ is the $S$-complex defined by
    \[
        \scC_\Lambda(S)
            =
        \bigcap \setst{ \cK \in \Com{S} }{ \textnormal{if $\cX \in \Lambda_\cK$, then $\bigcup \cX \in \cK$} };
    \]
that is, $\scC_\Lambda(S)$ is the minimal $S$-complex closed under unioning subsets assigned to it by the modulus $\Lambda$.
Equipping the object map $\scC_\Lambda$ with the action on morphisms sending $\varphi : S \rightarrow T$ to the extension
    \[
        \ext{\varphi} : \scC_\Lambda(S) \arw \scC_\Lambda(T)
            \quad \text{given by} \quad
        (X)\ext{\varphi} = \setst{(x)\varphi}{x \in X}
    \]
yields a functor which, moreover, admits a monad structure $(\scC_\Lambda , \sigma_\Lambda , \mu_\Lambda)$,
where the components of the unit $\sigma_\Lambda : \id{\FSGP} \Rightarrow \scC_\Lambda$ are the singleton embeddings
    \[
        \sigma_{\Lambda , S} 
            =
        \{-\} : S \longinj \scC_\Lambda(S)
            \quad \text{given by} \quad
        x \longmapsto \{x\}
    \]
and the components of the multiplication $\mu_\Lambda : \scC_\Lambda^2 \Rightarrow \scC_\Lambda$ are the union maps
    \[
        \mu_{\Lambda , S} 
            =
        \bigcup (-) : \scC_\Lambda^2(S) \longsurj \scC_\Lambda(S)
            \quad \text{given by} \quad
        \cX \longmapsto \bigcup_{X \in \cX} X
    \]
for every finite semigroup $S$.

The set of \textbf{points} of a modulus $\Lambda$, which is defined by
    \[
        \points{\Lambda}    
            \;=\;
        \setst{S \in \FSGP}{\Lambda_S \subseteq \sing(S)},
    \]
is a pseudovariety (see \cite[Proposition~9.7]{GENTHEORYPL}) with the additional property that $S \in \points{\Lambda}$ if and only if $\scC_\Lambda(S) = \sing(S)$.

This concept's utility comes from \cite[Theorem~9.12]{GENTHEORYPL}, which states that if $\Lambda$ is a modulus with $\points{\Lambda} = \pvV$, then $\scC_\Lambda(S) \subseteq \PV(S)$ for all $S \in \FSGP$.

\begin{notation}
If $\Lambda$ is a modulus and $\pvV$ is a pseudovariety, we write $\scC_\Lambda \leq \PV$ to indicate that $\scC_\Lambda(S) \subseteq \PV(S)$ for all finite semigroups $S$, and we write $\PV \leq \scC_\Lambda$ to mean the evident analogous statement.
\end{notation}

%
%%
%%%
%%%%

%%%%
%%%
%%
%

\section{Lower bound}   \label{section:lowerbound}

\subsection{Modulus}
Define a modulus $\Lambda_\pvER$ by
    \[
        \Lambda_\pvER = \modrl{S}{ \setst{ \BII{e} }{ e \in E(S) } },
    \]
and let $\CER(S)$ denote the $\Lambda_\pvER$-construct (\ref{framework:lowerbound}) of a given finite semigroup $S$.

\begin{lem}
The rule $\Lambda_\pvER$ is a modulus.
\end{lem}

\begin{proof}
Let $S$ be a finite semigroup.
Lemma \ref{type2:lem:type2KG} implies that $\BII{e} = \KG(S) \cap \RCL{e}$ for any $e \in E(S)$, from which it follows by basic stuff that $\BII{e}$ is a regular $\GR$-class of $\KG(S)$.
Since $\KG$ is an endofunctor which preserves regular epimorphisms by Lemma \ref{type2:lem:kgepi}, the required axioms are easily verified.
\end{proof}

\begin{prop}    \label{prop:lowerbound}
There is an equality $\pvER = \points{\Lambda_\pvER}$, and thus $\CER \leq \PER$.
\end{prop}

\begin{proof}
Let $S$ be a finite semigroup.
By Lemma \ref{type2:lem:ERchar}, $S \in \pvER$ if and only if the PTS $(R, S)$ is injective for every $\GR$-class $R$ of $S$.
Since $\equiv_\II$ is the minimal injective congruence on each $(R, S)$ by Lemma \ref{type2:lem:mininj}, this condition is equivalent to the condition that $\equiv_\II$ is the identity relation on all of $S$.

It is clear that if $\BII{x} = \{x\}$ for all $x \in S$ then $S \in \points{\Lambda_\pvER}$, from which it follows that $\pvER \subseteq \points{\Lambda_\pvER}$.

For the converse, suppose that $S \in \points{\Lambda_\pvER}$.
Given $x \in S$, Lemma \ref{type2:lem:actII} states that there exists $e \in E(S)$ such that $x \cdot \BII{e} = \BII{x}$.
But ex hypothesi $\BII{e} = \{e\}$, and hence $\BII{x} = \{x\}$ as well.
Thus $\points{\Lambda_\pvER} \subseteq \pvER$ and the desired equality holds.

The conclusion then follows from \cite[Theorem~9.12]{GENTHEORYPL} (as discussed in \ref{framework:lowerbound}).
\end{proof}

%
%%
%%%
%%%%

%%%%
%%%
%%
%

\section{Type-II blocks of pointlike sets}  \label{section:stable}

\subsection{}
For the remainder of the paper, fix a finite semigroup $S$.

\begin{defn}
Define a map $\beta : \CER(S) \rightarrow \CER(S)$ by setting
    \[
        (X)\beta = \bigcup \BII{X}
    \]
for each $X \in \CER(S)$.
\end{defn}

\begin{lem}
The semigroup $\CER(S)$ is closed under the action of $\beta$.
\end{lem}

\begin{proof}
Let $X \in \CER(S)$.
By Lemma \ref{type2:lem:actII} there is an idempotent $E \in \CER(S)$ such that $\BII{X} = X \cdot \BII{E}$.
But $\bigcup \BII{E} = (E)\beta$ is an element of $\CER(S)$ by the definition of $\Lambda_\pvER$, and thus $(X)\beta = X \cdot (E)\beta$ belongs to $\CER(S)$ as well.
\end{proof}

\begin{lem}     \label{stable:lem:blowup}
If $X \in \CER(S)$, then $X \subseteq (X)\beta$ and $(X)\beta \lgR X$.
In particular, if $A \in \cF_X$ then $(X)\beta = X \cdot (A)\beta$.
\end{lem}

\begin{proof}
It is obvious that $X \subseteq (X)\beta$.
As for the second claim, Lemmas \ref{prelim:lem:activator} and \ref{type2:lem:actII} guarantee the existence of some $A \in \cF_X$ such that $\BII{X} = X \cdot \BII{A}$, from which it follows that $(X)\beta = X \cdot (A)\beta$.
\end{proof}

\begin{lem} \label{stable:lem:idptblowup}
Let $E \in \CER(S)$ be an idempotent. Then
\begin{enumerate}
\item $(E)\beta$ is aperiodic, that is, $(E\beta)^{\omega + 1} = (E\beta)^\omega$; and
\item if $(E\beta)^2 \gH (E)\beta$ then $(E\beta)^2 = (E)\beta$.
\end{enumerate}
\end{lem}

\begin{proof}
\hfill
\begin{enumerate}
\item Note that $E \cdot (E)\beta = (E)\beta$ and $E \subseteq (E)\beta$.
Hence
    \[
        (E)\beta 
            =
        E \cdot (E)\beta
            \subseteq
        (E)\beta \cdot (E)\beta;
    \]
and, more generally, if $(E\beta)^{k-1} \subseteq (E\beta)^k$ then
    \[
        (E\beta)^k
            =
        (E\beta)^{k-1} \cdot (E)\beta
            \subseteq
        (E\beta)^k \cdot (E)\beta
            =
        (E\beta)^{k+1}.
    \]
Hence, given a number $p$ such that $(E\beta)^{\omega + p} = (E\beta)^\omega$, one has that
    \[
        (E\beta)^\omega
            \subseteq
        (E\beta)^{\omega + 1}
            \subseteq
            \cdots
            \subseteq
        (E\beta)^{\omega + p - 1}
            \subseteq
        (E\beta)^{\omega + p} = (E\beta)^\omega,
    \]
from which we conclude that $(E\beta)^{\omega + 1} = (E\beta)^\omega$.
\item Since $(E)\beta$ is aperiodic, the $\GH$-class of $(E\beta)^\omega$ is a single point.
Hence if $(E\beta)^2 \gH (E)\beta$, then $(E)\beta \gH (E\beta)^\omega$ and so $(E\beta)^2 = (E)\beta$.
\end{enumerate}
\end{proof}

\subsection{}
For each $X \in \CER(S)$, choose an idempotent $E_X \in \cF_X$ such that if $F$ is an idempotent then $E_F = F$.

\begin{defn}
Define a map
    \[
        \psi : \CER(S) \arw \CER(S)
            \quad \text{given by} \quad
        (X)\psi = X \cdot (E_X \beta)^\omega.
    \]
Note that if $E$ is an idempotent then $(E)\psi = (E\beta)^\omega$.
\end{defn}

\begin{lem}     \label{stable:lem:psifacts}
Let $X \in \CER(S)$. Then
\begin{enumerate}
\item $(X)\psi \lgR X$;
\item $(X)\psi = (X)\beta \cdot (E_X)\psi$ and $X \subseteq (X)\psi$;
\item $\psi$ is aperiodic, that is, $\psi^{\omega + 1} = \psi^\omega$;
\item if $X \gR (X)\psi$ then $X \equiv_\II (X)\psi$; and
\item $(X)\beta = X$ if and only if $(X)\psi = X$.
\end{enumerate} 
\end{lem}

\begin{proof}
\hfill
\begin{enumerate}
\item Obvious.
\item Since $(E_X \beta)^{\omega + 1} = (E_X \beta)^\omega$ by Lemma \ref{stable:lem:idptblowup}, 
    \[
        (X)\psi
            \;=\;
        X \cdot (E_X \beta)^\omega
            \;=\;
        X \cdot (E_X \beta) \cdot (E_X \beta)^\omega
            \;=\;
        (X)\beta \cdot (E_X \beta)^\omega.
    \]
It follows that 
    \[
        X 
            \;=\;
        X \cdot E_X
            \;\subseteq\;
        (X)\beta \cdot (E_X \beta)^\omega
            \;=\;
        (X)\psi
    \]
since $X \subseteq (X)\beta$ and $E_X \subseteq (E_X \beta)^\omega$.
\item Straightforward.
\item Since $(E_X)\psi = (E_X \beta)^\omega$ is an idempotent, it belongs to $\KG(\CER(S))$.
Hence if $(X)\psi \gR X$ then 
    \[
        \BII{(X)\psi} 
            =
        \BII{X \cdot (E_X)\psi}
            =
        \BII{X} \ast (E_X)\psi
            =
        \BII{X}
    \]
since members of $\KG(\CER(S))$ act as identity on $\II$-classes when defined.
\item For the ``only if'' direction, observe that
    \[
        (X)\beta = X \cdot (E_X)\beta = X
            \quad \Longrightarrow \quad
        X \cdot (E_X \beta)^k = X
    \]
for all $k \geq 1$, and hence $(X)\psi = X \cdot (E_X \beta)^\omega = X$.

As for the ``if'' direction, it follows from Lemma \ref{prelim:lem:activator} that if $(X)\psi = X$ then $(E_X)\psi \gR E_X$.
Consequently,
    \[
        (E_X)\psi
            =
        (E_X \beta)^\omega
            =
        (E_X \beta)^2
            =
        (E_X)\beta
    \]
by way of Lemma \ref{stable:lem:idptblowup}.
Therefore
    \[
        (X)\beta
            =
        X \cdot (E_X)\beta
            =
        X \cdot (E_X)\psi
            =
        (X)\psi
            =
        X,
    \]
at which point all desired claims have been established.
\end{enumerate}
\end{proof}

\begin{defn}
For each $X \in \CER(S)$, let $\clos{X} = (X)\psi^\omega$.
\end{defn}

\subsection{Fixed point sets}
Let $\bbF$ denote the set of fixed points of $\psi$; that is, let
    \[
        \bbF 
            \;=\; 
        \setst{ X \in \CER(S) }{ (X)\psi = X}
            \;=\;
        \setst{ \clos{X} }{ X \in \CER(S) },
    \]
which is also the set of fixed points of $\beta$.
Moreover, define
    \[
        \bbB 
            \;=\;
        \setst{ \BII{X} }{ X \in \bbF }
            \;=\;
        \setst{ \pi \in \CER(S)/\II }{ \bigcup \pi \in \pi }.
    \]

\begin{lem}     \label{stable:lem:stabilitylem}
Let $X \in \bbF$.
If $Y \in \CER(S)$ for which $(XY)\psi \gR XY \gR X$, then it follows that $\BII{XY} = \BII{X} \ast Y \in \bbB$.
\end{lem}

\begin{proof}
It follows from claims (2) and (4) of Lemma \ref{stable:lem:psifacts} that
    \[
        (XY)\psi \gR XY \gR X
            \quad \Longrightarrow \quad
        XY \equiv_\II (XY)\psi \equiv_\II (XY)\beta
    \]
since $(XY)\psi = (XY)\beta \cdot (E_{XY} \beta)^\omega$ (and since idempotents act as partial identity on $\II$-blocks), and so $\BII{XY} = \BII{X} \ast Y \in \bbB$.
\end{proof}

%
%%
%%%
%%%%

%%%%
%%%
%%
%

\section{Automaton and flow}   \label{section:automaton}

\subsection{}
An \textbf{automaton} is the data of a tuple $\cA = (\Sigma, Q, \INIT, \INP)$, where
\begin{enumerate}
    \item $\Sigma$ is a finite set of \textbf{input symbols},
    \item $Q$ is a finite set of \textbf{states},
    \item $\INIT \in Q$ is a distinguished \textbf{initial state}, and
    \item $(- \INP -) : Q \times \Sigma \rightarrow Q$ is a set function called the \textbf{transition function}.
\end{enumerate}
The \textbf{transition semigroup} of $\cA$ is the semigroup $\scT_\cA$ which is generated by the functions $(- \INP a) :  Q \rightarrow Q$ induced by each $a \in \Sigma$.

\begin{defn}
A \textbf{flow automaton} is a triple $(S, \cA, \Phi)$ where $S$ is a finite semigroup, $\cA = (S, Q, \INIT, \INP)$ is an automaton, and $\Phi$ (the nominal \textit{flow}) is a set function
    \[
        \Phi : Q \setminus \{ \INIT \} \arw \Po(S),
    \]
such that
    \[
    s \in (\INIT \INP s)\Phi
    \qquad \textnormal{and} \qquad
    (q)\Phi \cdot \{s\} \;\subseteq\; (q \INP s)\Phi
    \]
for all $s \in S$ and all $q \in Q \setminus \{\INIT\}$.
\end{defn}

\subsection{}       \label{automaton:covercomplex}
The \textbf{cover complex} of a flow automaton $(S, \cA, \Phi)$ is defined by
    \[
        \Cov(S , \cA, \Phi) 
            = 
        \setst{X \in \Po(S)}{\text{$X \subseteq (q)\Phi$ for some $q \in Q \setminus \{\INIT\}$}}.
    \]
It is straightforward to see that
    \[
        \PV(S)
        \;=\;
        \bigcap \setst{\Cov(S , \cA, \Phi)}{\scT_\cA \in \pvV}
    \]
for any pseudovariety $\pvV$ (see \cite[Proposition~2.5]{PL-VARDETGRP-BEN}).

\subsection{}
We will prove that $\PV(S) \subseteq \CER(S)$ by defining a flow automaton
    \[
        (S , \cA(S), \Phi)
            \quad \text{where} \quad
        \cA(S) = (S, Q(S), \INIT, \INP)
    \]
such that $\Cov(S , \cA(S), \Phi) \subseteq \CER(S)$ and $\scT_{\cA(S)} \in \pvER$.

\subsection{Local group actions}
Let $R$ be an $\GR$-class of $\CER(S)$, and consider the PTS $(R/\II , S)$ (which is isomorphic to $(R/\II , \sing(S))$ in the evident manner).
For each $s \in S$, extend the partial injection $(-) \ast s$ to a permutation $g_{(s , R)}$ on $R/\II$, whose action is written as
    \[
        \BII{X} \xmapsto{\;\;g_{(s, R)}\;\;} \BII{X} \GACT g_{(s, R)}
    \]
for each $\BII{X} \in R/\II$, and which has the property that
    \[
        \BII{X} \ast s \neq \varnothing
            \quad \Longrightarrow \quad
        \BII{X} \GACT g_{(s, R)} = \BII{X} \ast s
    \]
for all $\BII{X} \in R/\II$.
Let $G_R$ denote the group of permutations of $R/\II$ generated by the various $g_{(s ,R)}$ as $s$ ranges over $S$.

\subsection{Global group actions}
Let $\scR$ denote the set of $\GR$-classes of $\CER(S)$ which contain some member of $\bbF$; that is,
    \[
        \scR 
            =
        \setst{ \RCL{X} }{ (X)\psi = X }.
    \]
For each $s \in S$, let $g_s = (g_{(s , R)})_{R \in \scR}$.
Moreover, let
    \[
        \bbG = \left\langle g_s \mid s \in S \right\rangle;
    \]
that is, $\bbG$ is the subsemigroup of $\prod_{R \in \scR} G_R$ generated by the various $\scR$-tuples $g_s$ as $s$ ranges over $S$.
The group $\bbG$ acts on $\bigcup_{R \in \scR} R / \II$ by
    \[
        \BII{X} \GACT g 
            =
        \BII{X} \GACT g_{\RCL{X}}
    \]
for all $\BII{X} \in \bigcup_{R \in \scR} R / \II$ and all $g = (g_R)_{R \in \scR} \in \bbG$.

\subsection{States}
The state set of our automaton will be given by
    \[
        Q(S)
            =
        \setst{ (X , d , g) \in \bbF \times \bbG \times \bbG }{ \BII{X} \GACT d^\inv g \in \bbB }
            \cup
        \{ \INIT \}.
    \]
Moreover, define a map
    \[
        \val{-} : Q(S) \setminus \{ \INIT \} \arw \CER(S)
            \quad \text{given by} \quad
        \val{X , d , g} 
            =
        \bigcup \left( \BII{X} \GACT d^\inv g \right).
    \]
Note that the definition of $Q(S)$ guarantees that $\val{X , d , g} \in \bbF$ always.

\subsection{Updating maps}
For each $s \in S$, define a map $\lambda_s : Q(S) \setminus \{\INIT\} \rightarrow \bbF \times \bbG$ by 
    \[
        (X , d , g)\lambda_s 
            =
        \begin{cases}
            (X , d),
                    &   \text{if $\val{X , d , g} \cdot \{s\} \gR \val{X , d , g}$ and $(X , d , gg_s) \in Q(S)$;}
                    \\
            \left( \clos{ \val{X , d , g} \cdot \{s\} } , gg_s \right),
                    &   \text{otherwise;}
        \end{cases}
    \]
for each non-initial state $(X , d , g) \in Q(S) \setminus \{\INIT\}$.

\subsection{Automaton action}
The action of $s \in S$ in our automaton is given by
    \[
        \INIT \INP s 
            = 
        \left( \clos{ \{s\} } , g_s , g_s \right)
            \quad \text{and} \quad
        (X , d , g) \INP s
            =
        ( (X , d , g)\lambda_s , gg_s )
    \]
at each non-initial state $(X , d , g) \in Q(S) \setminus \{\INIT\}$.

\begin{lem}     \label{automaton:lem:actionwelldef}
If $q \in Q(S)$ and $s \in S$ then $q \INP s \in Q(S) \setminus \{ \INIT \}$.
\end{lem}

\begin{proof}
Clearly $q \INP s \neq \INIT$ always.
Observe that $q \INP s$ is always of the form
    \[
        q \INP s 
            =
        (X , d , gg_s)
    \]
for some $X \in \bbF$ and some $d , g \in \bbG$.
Notice that there are two possible cases.
\begin{enumerate}
\item In one case it is guaranteed that $q \INP s = (X , d, gg_s) \in Q(S)$.
\item Otherwise one has that $d = gg_s$, and so
    \[
        \val{ q \INP s } 
            =
        \bigcup \left( \BII{X} \GACT (gg_s)^\inv (gg_s) \right)
            =
        \bigcup \BII{X}
            =
        X.
    \]
It follows immediately that $q \INP s \in Q(S)$.
\end{enumerate}
Since these are the only two possibilities, we are done.
\end{proof}

\subsection{Flow automaton}
The data of our automaton $\cA(S)$ is given by
    \[
        \cA(S) = (S, Q(S), \INIT, \INP)
    \]
as defined thusfar in the section.
Our flow $\Phi$ will be the map 
    \[
        \Phi = \val{-} : Q(S) \setminus \{ \INIT \} \arw \Po(S),
    \]
whose image is clearly contained in $\CER(S)$.

\begin{prop}        \label{automaton:prop:flow}
The map $\Phi$ is a flow, and so $(S , \cA(S) , \Phi)$ is a flow automaton for which $\Cov(S , \cA(S), \Phi) \subseteq \CER(S)$.
\end{prop}

\begin{proof}
Given $s \in S$, we will consider $\val{ q \INP s }$ for each $q \in Q(S)$.
In the case where $q = \INIT$, the value of $\val{ \INIT \INP s }$ is equal to
    \[
        \val{ \left( \clos{ \{s\} } , g_s , g_s \right) }
            =
        \bigcup \left( \BII{ \clos{ \{s\} } } \GACT g_s^\inv g_s \right)
            =
        \bigcup \BII{ \clos{ \{s\} } }
            =
        \clos{ \{s\} },
    \]
of which $s$ is clearly an element.

The remaining two cases involve non-initial states $(X , d , g) \in Q(S) \setminus \{ \INIT \}$.
\begin{enumerate}
\item If $\val{X , d , g} \cdot \{s\} \gR \val{X , d , g}$ and $(X , d , gg_s) \in Q(S)$, then
    \[
        (X , d , g) \INP s = (X , d , gg_s)
            \quad \text{and} \quad
        \BII{X} \GACT d^\inv g g_s = \BII{ \val{X , d , g} \cdot \{s\} }.
    \]
From here, the computation
    \[
        \val{ (X , d , g) \INP s }
            =
        \bigcup \left( \BII{X} \GACT d^\inv g g_s \right)
            =
        \bigcup \BII{ \val{X , d , g} \cdot \{s\} }
    \]
yields the fact that $\val{X , d , g} \cdot \{s\} \subseteq \val{ (X , d , g) \INP s }$.
\item In all remaining cases one has that 
    \[
        (X , d , g) \INP s = \left( \clos{ \val{X , d , g} \cdot \{s\} } , gg_s , gg_s \right).
    \]
It follows from claim (2) of Lemma \ref{stable:lem:psifacts} that
    \[
        \val{ (X , d , g) \INP s }
            =
        \clos{ \val{X , d , g} \cdot \{s\} }
            \supseteq
        \val{X , d , g} \cdot \{s\}.
    \]
\end{enumerate}
Having covered all cases, we have established that $\Phi$ is a flow.
\end{proof}

%
%%
%%%
%%%%

%%%%
%%%
%%
%

\section{Upper bound}   \label{section:upperbound}

\subsection{}
Given a preordered set $P$, let $\mathcal{D}_P$ denote the set of (not necessarily monotone) functions $f : P \rightarrow P$ which satisfy
    \[
        (x)f = x
            \quad \text{or} \quad
        (x)f < x
    \]
for all $x \in P$.
Clearly $\mathcal{D}_P$ is closed under composition and thus is a semigroup.

\begin{lem}     \label{upperbound:lem:rdown}
If $P$ is a preordered set, then $\mathcal{D}_P$ is $\GR$-trivial.
\end{lem}

\begin{proof}
Let $f, g, h \in \mathcal{D}_P$ and suppose that $fgh=f$.
If $x \in P$, then 
    \[
        (x)f = (x)fgh \leq (x)fg \leq (x)f,
    \]
and hence $(x)fg = (x)f$ already.
\end{proof}

\begin{defn}
Let $(\bbF \times \bbG)^\bullet = (\bbF \times \bbG) \cup \{ \bullet \}$.
Define a preorder $\leq$ on $(\bbF \times \bbG)^\bullet$ by setting
\begin{enumerate}
\item $\bullet \leq \bullet$ and $(X , d) \leq \bullet$ always; and
\item $(X_1 , d_1) \leq (X_2 , d_1)$ if and only if $X_1 \lgR X_2$.
\end{enumerate}
\end{defn}

\subsection{}       \label{upperbound:extendpointer}
Given $s \in S$ and $g \in \bbG$, extend the map $( - , g )\lambda_s$ to $(\bbF \times \bbG)^\bullet$ by setting
    \[
        ( \bullet , g )\lambda_s 
            =
        ( \clos{ \{s\} } , g_s );
    \]
and, for any $(X , d)$ such that $(X , d , g)$ does not belong to $Q(S)$, setting
    \[
        (X , d , g)\lambda_s = (X , d).
    \]

\begin{lem}     \label{upperbound:lem:pointerRdown}
The map $( - , g)\lambda_s$ belongs to $\mathcal{D}_{(\bbF \times \bbG)^\bullet}$ for all $s \in S$ and all $g \in G$.
\end{lem}

\begin{proof}
We begin by considering the ``extended'' cases defined in \ref{upperbound:extendpointer}.
First, since $\bullet$ is strictly above all non-$\bullet$ elements of $(\bbF \times \bbG)^\bullet$, it follows that $(\bullet , g)\lambda_s < \bullet$.
Next, if $(X , d , g)$ is not a member of $Q(S)$, then $(X , d , g)\lambda_s = (X , d)$.

We now move on to the cases where $(X , d , g) \in Q(S)$, of which there are three.
\begin{enumerate}
\item If $\val{X , d , g} \cdot \{s\} \gR \val{X , d , g}$ and $(X , d , gg_s) \in Q(S)$, then $(X , d , g)\lambda_s = (X , d)$.
\item If $\val{X , d , g} \cdot \{s\} \slgR \val{X , d , g}$, then 
    \[
         \clos{ \val{X , d , g} \cdot \{s\} } 
            \lgR
         \val{X , d , g} \cdot \{s\}
            \slgR
        \val{X , d , g}
            \gR
        X
    \]
and therefore $(X , d , g)\lambda_s < (X , d)$.
\item Finally, suppose that $\val{X , d , g} \cdot \{s\} \gR \val{X , d , g}$ but $(X , d , gg_s)$ does not belong to $Q(S)$.
In this situation, $\val{X , d , g} \in \bbF$ but the $\II$-block
    \[
        \BII{X} \GACT d^\inv g g_s
            =
        \BII{ \val{X , d , g} } \GACT g_s
            =
        \BII{ \val{X , d , g} \cdot \{s\} }
    \]
does not belong to $\bbB$.
This implies via the contrapositive of Lemma \ref{stable:lem:stabilitylem} that
    \[
        \clos{ \val{X , d , g} \cdot \{s\} }
            \slgR
        \val{X , d , g} \cdot \{s\}
            \gR
        \val{X , d , g}
            \gR 
        X,
    \]
from which we conclude that $(X , d , g)\lambda_s < (X , d)$.
\end{enumerate}
Having verified the desired condition in all cases, we are done.
\end{proof}

\begin{prop}    \label{upperbound:prop:computingER}
The transition semigroup of $\cA(S)$ belongs to $\pvER$.
\end{prop}

\begin{proof}
Let $\scT(S)$ denote the transition semigroup of $\cA(S)$.
Moreover, for each $s \in S$, let $\tilde{s}$ denote the transformation $(-) \INP s \in \scT(S)$, and note that $\scT(S)$ is generated by the various $\tilde{s}$ as $s$ ranges over $S$.

We will define an embedding of (faithful) transformation semigroups
    \[
        (\zeta , \varphi) : (Q(S) , \scT(S)) 
            \longinj 
        ( (\bbF \times \bbG)^\bullet , \mathcal{D}_{(\bbF \times \bbG)^\bullet})
            \wr
        (\bbG , \bbG),
    \]
which will establish via Lemmas \ref{upperbound:lem:rdown} and \ref{type2:lem:ERchar} that $\scT(S) \in \pvR \ast \pvG = \pvER$.

The map $\zeta : Q(S) \rightarrow (\bbF \times \bbG)^\bullet \times \bbG$ is given by
    \[
        (\INIT)\zeta = (\bullet , 1_{\bbG})
            \quad \text{and} \quad
        (X , d , g)\zeta = ( X , d , g).
    \]
It is clear that $\zeta$ is injective.

Next, we define the morphism $\varphi : \scT(S) \rightarrow \mathcal{D}_{(\bbF \times \bbG)^\bullet} \wr (\bbG , \bbG)$.
To do so, it suffices to define $(\tilde{s})\varphi$ for all $s \in S$ since $\varphi$ is determined by its values on the generators of $\scT(S)$.
So, for each $s \in S$, define
    \[
        (\tilde{s})\varphi = (\tilde{\lambda}_s , g_s)
    \]
where the function $\tilde{\lambda}_s : \bbG \rightarrow \mathcal{D}_{(\bbF \times \bbG)^\bullet}$ is given by
    \[
        (g)\tilde{\lambda}_s 
            = 
        ( - , g)\lambda_s : (\bbF \times \bbG)^\bullet \rightarrow (\bbF \times \bbG)^\bullet
    \]
at each group element $g \in \bbG$.
This is well-defined since each $( - , g)\lambda_s$ belongs to $\mathcal{D}_{(\bbF \times \bbG)^\bullet}$ by Lemma \ref{upperbound:lem:pointerRdown}.
Moreover, it is clear that $\varphi$ is injective.

Now, observe that
    \[
        (\INIT)\zeta \ast (\tilde{s})\varphi 
            \;=\;
        \left( \bullet , 1_{\bbG} \right) \ast (\tilde{\lambda}_s , g_s)
            \;=\;
        \left( \clos{ \{s\} } , g_s , g_s \right)
            \;=\;
        (\INIT \INP s)\zeta;
    \]
and if $(X , d, g) \in Q(S) \setminus \{ \INIT \}$ then
    \[
        (X , d, g)\zeta \ast (\tilde{s})\varphi 
            \;=\;
        ( (X , d , g)\lambda_s , g g_s )
            \;=\;
        ( (X , d , g) \INP s )\zeta.
    \]
This establishes that the pair $(\zeta , \varphi)$ is an embedding of transformation semigroups.

Since $\mathcal{D}_{(\bbF \times \bbG)^\bullet}$ is $\GR$-trivial by Lemma \ref{upperbound:lem:rdown}---and since $\bbG$ is obviously a group---it follows that $\scT(S)$ belongs to the pseudovariety $\pvR \ast \pvG$.
Since $\pvR \ast \pvG = \pvER$ by Lemma \ref{type2:lem:ERchar}, the proposition follows.
\end{proof}

\begin{thm}     \label{thm:maintheorem}
Pointlike sets are decidable for $\pvER$.
In particular, $\PER = \CER$.
\end{thm}

\begin{proof}
It was established in Proposition \ref{prop:lowerbound} that $\CER \leq \PER$.
As for the other bound, Propositions \ref{automaton:prop:flow} and \ref{upperbound:prop:computingER} show that for each finite semigroup $S$ there exists a flow automaton whose cover complex is contained in $\CER(S)$ and whose transition semigroup belongs to $\pvER$.
Considering this alongside \ref{automaton:covercomplex} establishes that $\PER \leq \CER$.
Thus $\PER = \CER$; and, since $\CER$ is computable, we conclude that $\pvER$ has decidable pointlikes.
\end{proof}

%
%%
%%%
%%%%

\printbibliography
\hrulefill
\end{document}